\documentclass[psamsfonts,reqno]{amsart}
\usepackage{amssymb,eucal,latexsym}
\newcommand{\MM}{\mathcal{M}}

\newcommand{\les}{\leqslant}
\newcommand{\ges}{\geqslant}
\newcommand{\nor}[1]{{#1}^{\sharp}}
\newcommand{\ol}[1]{\overline{#1}}
\newcommand{\supp}{\mathrm{supp}}
\newcommand{\bsupp}{\ol{\supp}}
\newcommand{\dnw}{\mathbin{\downarrow}}

\newcommand{\pI}[1]{\bigl(#1\bigr)}

\newcommand{\del}[2]{\Delta_{#1}(#2)}

\newcommand{\jirr}{join-ir\-re\-duc\-i\-ble}

\newcommand{\jsd}{join-sem\-i\-dis\-trib\-u\-tive}

\newcommand{\contr}{a contradiction}

\newcommand{\DD}{\mathbin{D}}

\newcommand{\pup}[1]{\textup{(}{#1}\textup{)}}

\newcommand{\fsi}[1]{\{1,\dots,#1\}}
\newcommand{\fso}[1]{\{0,\dots,#1\}}

\newcommand{\es}{\varnothing}

\newcommand{\ad}{{\dot{a}}}
\newcommand{\bd}{{\dot{b}}}
\newcommand{\cd}{{\dot{c}}}
\newcommand{\pd}{{\dot{p}}}
\newcommand{\qd}{{\dot{q}}}
\newcommand{\rd}{{\dot{r}}}
\newcommand{\xd}{{\dot{x}}}

\newcommand{\tr}{\vartriangleleft}

\newcommand{\utr}{\trianglelefteq}

\newcommand{\dtr}{\mathbin{\vartriangleleft\kern-10pt
{\lower3pt\hbox{$\scriptscriptstyle\neq$}}\kern3pt}}

\newcommand{\eb}{\partial_{\mathrm{e}}}
\newcommand{\into}{\hookrightarrow}
\newcommand{\onto}{\twoheadrightarrow}
\newcommand{\twohead}{\mathbin{{-}\kern-2pt{\onto}}}

\newcommand{\lar}[1]{{\longrightarrow^{{#1}}}}
\newcommand{\larr}[1]{{\twohead^{{#1}}\kern3pt}}

\newcommand{\arr}[1]{\mathbin{\underset{(#1)}{\longrightarrow}}}
\newcommand{\arrf}[1]{\mathbin{\underset{(#1)}{\twohead}}}

\newcommand{\set}[1]{\{{#1}\}}
\newcommand{\setm}[2]{\set{{#1}\mid{#2}}}
\newcommand{\seq}[1]{\langle#1\rangle}
\newcommand{\Set}[1]{\left\{{#1}\right\}}
\newcommand{\Setm}[2]{\Set{{#1}\mid{#2}}}

\newcommand{\two}{\mathbf{2}}
\newcommand{\ZZ}{\mathbb{Z}}
\newcommand{\QQ}{\mathbb{Q}}
\newcommand{\RR}{\mathbb{R}}
\newcommand{\FF}{\mathbb{F}}
\newcommand{\VV}{\mathbb{V}}
\newcommand{\eps}{\varepsilon}

\newcommand{\FFT}{\FF^{(T)}}
\newcommand{\FFTp}{\FF_+^{(T)}}

\DeclareMathOperator{\J}{J}

\newcommand{\hgt}{\mathrm{ht}}
\newcommand{\Co}{\mathbf{Co}}
\DeclareMathOperator{\Col}{Co}

\newcommand{\Case}[1]{\smallskip\noindent{\textbf{Case~{#1}}}}

\numberwithin{equation}{section}

\theoremstyle{plain}
\newtheorem{lemma}{Lemma}[section]
\newtheorem{theorem}[lemma]{Theorem}
\newtheorem{proposition}[lemma]{Proposition}
\newtheorem{corollary}[lemma]{Corollary}
\newtheorem{claim}{Claim}

\newtheorem*{stat}{\name}
\newcommand{\name}{testing}

\theoremstyle{definition}
\newtheorem{definition}[lemma]{Definition}
\newtheorem{notation}[lemma]{Notation}
\newtheorem{example}[lemma]{Example}
\newtheorem{problem}{Problem}

\theoremstyle{remark}
\newtheorem{remark}[lemma]{Remark}

\newenvironment{scproof}
{\begin{proof}[Proof of Claim.]}
{\end{proof}}

\begin{document}

\title[Lattices of convex sets]
{Sublattices of lattices of convex subsets of vector spaces}

\author[F.~Wehrung]{Friedrich Wehrung}
\address[F.~Wehrung]{CNRS, UMR 6139\\
D\'epartement de Math\'ematiques\\
Universit\'e de Caen\\
14032 Caen Cedex\\
France}
\email{wehrung@math.unicaen.fr\\
http://www.math.unicaen.fr/\~{}wehrung}

\author[M.\,V.~Semenova]{Marina Semenova}
\address[M.\,V.~Semenova]{Institute of Mathematics of
the Siberian Branch of RAS\\
Acad. Koptyug prosp. 4\\
630090 Novosibirsk\\
Russia}
\email{semenova@math.nsc.ru}

\date{\today}

\thanks{The authors were partially supported by GA CR grant
201/00/0766 and by institutional grant MSM:J13/98:1132000007a.
The second author was partially supported by INTAS grant
YSF: 2001/1-65, by the joint RFBR--DFG grant 01-01-04003 NNIOa,
by the Russian Ministry of Education grant E 02-1.0-32, by
the President of Russian Federation grant NSh-2112.2003.1 supporting
leading scientific schools, and by the Science Support Foundation grant.}

\begin{abstract}
For a left vector space $V$ over a totally ordered division ring $\FF$,
let $\Co(V)$ denote the lattice of convex subsets of $V$. We prove that
every lattice $L$ can be embedded into $\Co(V)$ for some left
$\FF$-vector space~$V$. Furthermore, if $L$ is finite lower bounded, then $V$
can be taken finite-dimensional, and $L$ embeds into a finite lower bounded
lattice of the form $\Co(V,\Omega)=\setm{X\cap\Omega}{X\in\Co(V)}$,
for some finite subset $\Omega$ of $V$. In particular, we obtain a
new universal class for finite lower bounded lattices.
\end{abstract}

\begin{flushright}
UDC 512.56
\end{flushright}

\maketitle

\section{Introduction}\label{S:Intro}
The question about the possibility to embed lattices from a particular class
into lattices from another particular class (or, the question about description
of sublattices of lattices from a particular class) has a long history.
Many remarkable
results were obtained in that direction. Among the first classical
ones, one can
mention the result of Ph.\,M. Whitman~\cite{W} published in 1946 that
every lattice
embeds into the partition lattice of a set. The question whether every
\emph{finite} lattice embeds into the partition lattice of a \emph{finite}
set was a long-standing problem, which was solved in the positive in 1980
by P. Pudl\'{a}k and J. T\r uma in their well-known paper~\cite{PuTu}.

The paper~\cite{AGT} by K.\,V. Adaricheva, V.\,A. Gorbunov, and V.\,I.
Tumanov investigates the question of embedding lattices into so-called
\emph{convex geometries}, that is, closure lattices of closure spaces
with the anti-exchange property.
It is well-known that any finite convex geometry is \emph{\jsd}, that is, it
satisfies the following quasi-identity:
       \[
       \forall xyz\ x\vee y=x\vee z\rightarrow x\vee y=x\vee(y\wedge z).
       \]
Moreover, it is proved in~\cite[Theorem~1.11]{AGT} that
any finite \jsd\ lattice embeds into a finite convex geometry. Among
other things,
one particular class of convex geometries, the class of lattices of algebraic
subsets of complete lattices, was studied in the abovementioned paper.
The authors of~\cite{AGT} proved that any finite \jsd\ lattice embeds into
the lattice of algebraic subsets of some algebraic and dually algebraic
complete lattice~$A$. In general, the lattice~$A$ may be infinite.
This result inspired Problem~3 in~\cite{AGT}, which asks the following:
    \begin{quote}\em
    Is there a special class $\mathcal{U}$ of finite convex geometries
    that contains all finite \jsd\ lattices as sublattices?
    \end{quote}
In other words, is there a special class $\mathcal{U}$ of finite convex
geometries such that any finite \jsd\ lattice embeds into a lattice from
$\mathcal{U}$?
For the class of subsemilattice lattices of finite semilattices,
an answer to the above
question is provided by the following result which was proved
independently by K.\,V. Adaricheva~\cite{Adar96} and
V.\,B. Repnitskii~\cite{R}:
    \begin{quote}\em
    A finite lattice embeds into the
    subsemilattice lattice of a finite \pup{semi}lattice iff it is
    lower bounded.
    \end{quote}
Another result of the same spirit was proved by B. \v{S}ivak~\cite{S}
(see also~\cite{Se}):
    \begin{quote}\em
    A finite lattice embeds into the
    suborder lattice of a finite partially ordered set iff it is
    lower bounded.
    \end{quote}
We observe that the class of finite lower bounded lattices is a proper
subclass of the class of finite \jsd\ lattices (see~\cite{FJN}).
For a precise definition of a lower bounded lattice, we refer the reader to
Section~\ref{S:Basic}.

As natural candidates for~$\mathcal{U}$, the following classes were proposed
in~\cite{AGT}:
\begin{itemize}
\item[(1)]
The class of all finite, atomistic, \jsd, biatomic lattices.
\item[(2)]
The class of all lattices of the form $\Co(P)$, the lattice of all
order-convex subsets of a finite partially ordered set $P$.
\item[(3)]
The class of all lattices of the form
$\Co(\RR^n,\Omega)=\setm{X\cap\Omega}{X\in\Co(\RR^n)}$,
for a finite $\Omega\subseteq\RR^n$ and $n<\omega$ (see
Section~\ref{S:Basic} for the notation).
\end{itemize}

The class (1) turns out to be too restrictive, see K.\,V. Adaricheva and F.
Weh\-rung~\cite{AdWe}.
The class (2) is even more restrictive. In~\cite{SeWe1}, the
sublattices of finite lattices of the form $\Co(P)$
are described; in particular, they are the finite lattices satisfying three
identities, denoted there by (S), (U), and (B).
Whether the class~(3) can be such a ``universal'' class $\mathcal{U}$
for finite \jsd\ lattices is still open (see Problem 1).

In the present paper, we prove that every lattice embeds into the
lattice of convex
subsets of a vector space (see Theorem~\ref{T:ArbLatt}).
We also get the following partial confirmation of the hypothesis about
``universality'' of the class (3) (see
Theorem~\ref{T:FinLB}):
    \begin{quote}\em
    Every finite lower bounded lattice embeds into $\Co(\RR^n,\Omega)$, for some
    $n<\omega$ and some finite $\Omega\subseteq\RR^n$.
    \end{quote}

Both main results of the paper are proved by using the same method of
construction, which is elaborated in Sections~\ref{S:ColTrees} to
\ref{S:Norms}. All the vector spaces that we shall consider will be built up
from so-called \emph{colored trees}, see Definition~\ref{D:ColTree}. The
elements of the tree have to be thought of as finite sequences of \jirr\
elements of the lattice we are starting from, together with some additional
information, as shown
in Section~\ref{S:EmbThms}. A precursor for this method
can be found in~\cite{Adar96}, where a meet-semilattice is constructed from
finite sequences of \jirr\ elements from the original finite lower bounded
lattice. See also Section~2.1 in~\cite{AGT}.
The elements of the tree $T$ index the canonical basis of the free
vector space $\FFT$ on $T$, and new relations on these elements are introduced
\emph{via} a \emph{rewriting rule}, denoted by $\lar{*}$, on the positive cone
$\FFTp$ of $\FFT$, see Section~\ref{S:ColTrees}. It turns out that this
rewriting rule is \emph{confluent} (Lemma~\ref{L:RevConfmn}), which makes it
possible to say that two elements are equivalent if{f} they have
some common rewriting, see Notation~\ref{No:DefEquiv} and
Proposition~\ref{P:EquivTrans}. Differences of equivalent elements form a
vector subspace, $N_T$, and the interesting convex subsets will live in the
vector space
$\VV_T=\FFT/{N_T}$, see
Section~\ref{S:CanThm}.

We offer two types of technical results.
Our first type of result states that equality of two elements of
$\FFT$ modulo $N_T$ can be conveniently expressed \emph{via} the rewriting
rule, essentially Proposition~\ref{P:EquivTrans} (the equivalence $\equiv$ can
be expressed \emph{via} common rewriting) and Theorem~\ref{T:EquivCan} (the
equivalence $\equiv$ is cancellative). These results are not lattice
theoretical, but combinatorial.

Our second type of result is more lattice theoretical, and it says
which sort of colored tree $T$ we need in order to embed a given lattice $L$
into $\Co(\VV_T)$ nicely. The most central result among those is
Theorem~\ref{T:Norm2Emb}. It uses the notion of a ``$L$-valued norm'' on a tree
$T$.

In Section~\ref{S:RelAbs}, we will show some relationship between embeddability
into $\Co(V)$ and into $\Co(V,\Omega)$. We conclude the paper with some open
problems in Section~\ref{S:Pbs}.

We observe
that the class of lattices of convex subsets of vector spaces
was studied by A. Huhn. In particular, he proved in~\cite{H} that, for a
$(n-1)$-dimensional vector space $V$, the lattice
$\Co(V)$ belongs to the variety generated by all finite
$n$-distributive lattices;
thus it is $n$-distributive itself, however, it is not $(n-1)$-distributive,
see also G.\,M. Bergman~\cite{Berg}.
In the finite dimensional case,
principal ideals of lattices of the form $\Co(V)$ are characterized in
M.\,K. Bennett~\cite{Benn77}.
An alternate proof of the second half of Theorem~\ref{T:FinLB}, that uses the
main result of~\cite{Adar96}, can be found in K.\,V. Adaricheva~\cite{CoX}.

\section{Basic concepts}\label{S:Basic}

We first recall some classical concepts,
about which we also refer the reader to  R.~Freese, J. Je\v{z}ek, and J.\,B.
Nation~\cite{FJN}. For a join-semilattice $L$,
we set $L^-=L\setminus\set{0}$ if $L$ has a zero
(least element), $L^-=L$ otherwise. For subsets $X$ and $Y$
of $L$, we write that $X\ll Y$, if every element of $X$ lies below some element
of~$Y$. If $a\in L^-$, a \emph{nontrivial join-cover} of $a$ is a
finite subset $X$ of $L^-$ such that $a\leq\bigvee X$ while $a\nleq x$ for all
$x\in X$. A nontrivial join-cover $X$ of $a$ is \emph{minimal}, if $Y\ll X$
implies that $X\subseteq Y$, for any nontrivial join-cover $Y$ of $a$.
We denote by $\J(L)$ the set of all \jirr\ elements of $L$.
For $a,b\in\J(L)$, we write $a\DD b$ if $b$ belongs to a minimal
nontrivial join-cover of $a$.
A sequence
$a_0,\ldots,a_{n-1}$ of elements from $\J(L)$ is a $\DD$-\emph{cycle}, if
$a_0\DD\ldots\DD a_{n-1}\DD a_0$.

A lattice homomorphism $h\colon K\to L$ is \emph{lower bounded}
if, for all $a\in L$, the set $\setm{x\in K}{h(x)\geq a}$ is
either empty or has a least element. A finitely generated lattice
$L$ is \emph{lower bounded}, if it is the homomorphic
image of a finitely generated
free lattice under a lower bounded lattice homomorphism.
Equivalently, for
finite~$L$, the $\DD$ relation of $L$ has no cycle.

For posets $K$ and $L$, we say that a map $f\colon K\to L$ is
\emph{zero-preserving}, if whenever~$K$ has a smallest element, say, $0_K$,
the element $f(0_K)$ is the smallest element of $L$. We say that $f$
\emph{preserves existing meets}, if whenever $X\subseteq K$ has a meet in
$K$, the image $f[X]$ has a meet in $L$, and
$\bigwedge f[X]=f\left(\bigwedge X\right)$.

For a totally ordered division ring $\FF$ and a positive integer $n$, we put
       \begin{equation}\label{Eq:DelnF}
       \del{n}{\FF}=\Setm{(\xi_i)_{i<n}\in(\FF^+)^n}{\sum_{i<n}\xi_i=1},
       \end{equation}
the \emph{$(n-1)$-simplex} in $\FF^n$.

All vector spaces considered in this paper will be left vector spaces.
Let $V$ be a vector space over a totally ordered division ring $\FF$. We put
       \[
       [x,y]=\setm{\xi_0x+\xi_1y}{(\xi_0,\xi_1)\in\del{2}{\FF}},
       \]
for all $x$, $y\in V$. A subset $X$ of $V$ is \emph{convex}, if
$[x,y]\subseteq X$ whenever $x$, $y\in X$. We denote by $\Co(V)$ the lattice
(under inclusion) of all convex subsets of~$V$.
For a subset $X$ of $V$, we denote by
$\Col(X)$ the \emph{convex hull} of~$X$. Hence
       \[
       \Col(X)=\Setm{\sum_{i<n}\xi_ix_i}
       {0<n<\omega,\ (\xi_i)_{i<n}\in\del{n}{\FF},\ (x_i)_{i<n}\in X^n}.
       \]
For a subset $\Omega$ of $V$, we put
       \[
       \Co(V,\Omega)=\setm{X\cap\Omega}{X\in\Co(V)}.
       \]
In general, $\Co(V,\Omega)$ is a lattice, it is, in fact, (the closure lattice
of) a \emph{convex geometry}, see~\cite{AGT}. As shows the following result,
there are only trivial \jirr\ elements in $\Co(V,\Omega)$, even for infinite
$\Omega$.

\begin{proposition}\label{P:JirrCoVO}
Let $V$ be a vector space over a totally ordered division ring $\FF$, let
$\Omega$ be a subset of $V$. Then the \jirr\ elements of $\Co(V,\Omega)$ are
exactly the singletons $\set{p}$, for $p\in\Omega$.
\end{proposition}

\begin{proof}
It is trivial that singletons of elements of $\Omega$ are (completely) \jirr.
Let $P$ be \jirr\ in $\Co(V,\Omega)$, suppose that there are distinct $a$,
$b\in P$. There exists a linear functional $f\colon V\to\FF$ such that
$f(a)<f(b)$. Put
      \begin{align*}
      X&=\setm{x\in P}{f(x)\les f(a)},\\
      Y&=\setm{x\in P}{f(x)>f(a)}.
      \end{align*}
Then $X$, $Y$ belong to $\Co(V,\Omega)$, $P=X\cup Y=X\vee Y$, and $X$,
$Y\neq P$, \contr.
\end{proof}

For a partially ordered abelian group $G$, we put
       \[
       G^+=\setm{x\in G}{0\leq x},\qquad G^{++}=G^+\setminus\set{0}.
       \]
We shall need a few elementary binary operations on \emph{ordinals}:
we denote by $(\alpha,\beta)\mapsto\alpha^{\beta}$ the exponentiation, by
$(\alpha,\beta)\mapsto\alpha\dotplus\beta$ the addition, by
$(\alpha,\beta)\mapsto\alpha\cdot\beta$ the multiplication, and
by $(\alpha,\beta)\mapsto\alpha+\beta$ the \emph{natural sum} (or
\emph{Hessenberg sum}), see K. Kuratowski and A. Mostowski~\cite{KuMo}.
By definition, if $k$, $n_0$, \dots, $n_{k-1}$,
$p_0$, \dots, $p_{k-1}$, $q_0$, \dots, $q_{k-1}$ are natural numbers such that
$n_0>n_1>\cdots>n_{k-1}$, and
       \begin{align*}
       \alpha&=\omega^{n_0}\cdot p_0\dotplus\cdots\dotplus
       \omega^{n_{k-1}}\cdot p_{k-1},\\
       \beta&=\omega^{n_0}\cdot q_0\dotplus\cdots\dotplus
       \omega^{n_{k-1}}\cdot q_{k-1},
       \end{align*}
then the Hessenberg sum of $\alpha$ and $\beta$ is given by
       \[
       \alpha+\beta=\omega^{n_0}\cdot(p_0+q_0)\dotplus\cdots\dotplus
       \omega^{n_{k-1}}\cdot(p_{k-1}+q_{k-1}).
       \]
In particular, the Hessenberg addition is commutative, associative,
and cancellative. Moreover,
if $n_0\geq n_1\geq\cdots\geq n_{k-1}$ are natural numbers, then the
Hessenberg sum of the $\omega^{n_i}$-s is given by
       \[
\sum_{i<k}\omega^{n_i}=\omega^{n_0}\dotplus\cdots\dotplus\omega^{n_{k-1}}.
       \]

\section{The free vector space associated with a colored tree}
\label{S:ColTrees}

Let $(T,\utr)$ be a partially ordered set. We denote by $\tr$ the
associated strict ordering of $T$. For elements $a$ and $b$ of $T$, we
say that $a$ is a \emph{lower cover} of $b$, in notation
$a\prec b$, if $a\tr b$ and there exists no element $x$ of $T$ such that
$a\tr x\tr b$. If~$b$ has exactly one lower cover, we denote it by $b_*$.

\begin{definition}\label{D:ColTree}
A \emph{tree} is a partially ordered set $(T,\utr)$ such that the lower
segment $\dnw p=\setm{q\in T}{q\utr p}$ is a finite chain, for any $p\in T$.
We put $\hgt(p)=|\dnw p|-1$, for all $p\in T$.

A \emph{coloring} of a tree $(T,\utr)$ is an equivalence relation $\sim$
on $T$ such that the following statements hold:
\begin{enumerate}
\item The $\sim$-equivalence class $[p]$ of $p$ is finite and has at
least two elements, for any non-minimal $p\in T$.

\item If $p\sim q$, then either both $p$ and $q$ are minimal or $p_*=q_*$, for
all $p$, $q\in T$.
\end{enumerate}
A \emph{colored tree} is a triple $(T,\utr,\sim)$, where $(T,\utr)$ is a
tree and $\sim$ is a coloring of $T$.
\end{definition}

For a colored tree $(T,\utr,\sim)$, we put
       \begin{align*}
       \MM_T&=\setm{(p,[q])}{p,\,q\in T\text{ and }p\prec q};\\
       \MM_T(p)&=\setm{[q]}{q\in T\text{ and }p\prec q},\text{ for all }p\in T.
       \end{align*}
For a totally ordered division ring $\FF$, we consider the free vector
space $\FFT$ on~$T$, whose elements are the maps $x\colon T\to\FF$
whose
\emph{support} $\supp(x)=\setm{p\in T}{x(p)\neq0}$ is finite. We denote by
$(\pd)_{p\in T}$ the canonical basis of $\FFT$, and we order $\FFT$
componentwise, that is,
       \[
       x\leq y,\text{ if }x(p)\les y(p)\text{ for all }p\in T.
       \]
With this ordering, $\FFT$ is a \emph{lattice-ordered vector space}
over $\FF$. We denote by $\FFTp$ the positive cone of $\FFT$, that is,
       \[
       \FFTp=\setm{x\in\FFT}{x(p)\ges 0\text{ for all }p\in T}.
       \]
For $(p,I)\in\MM_T$, we define binary relations $\arr{p,I}$ and
$\arrf{p,I}$ on $\FFTp$ by
       \begin{align*}
       &x\arr{p,I}y\ \Longleftrightarrow
       \text{ there are }\lambda\in\FF^+
       \text{ and }z\in\FFTp\text{ such that }\\
       &x=\frac{\lambda}{|I|}\sum_{q\in I}\qd+z\text{ and }y=\lambda\pd+z.\\
       &x\arrf{p,I}y\ \Longleftrightarrow
       \text{ there are }\lambda\in\FF^+
       \text{ and }z\in\FFTp\text{ such that }\\
       &x=\frac{\lambda}{|I|}\sum_{q\in I}\qd+z,\ y=\lambda\pd+z,\text{ and }
       z(q_0)=0\text{ for some }q_0\in I.
       \end{align*}
If $x\arr{p,I}y$, we say that $y$ is the result of a \emph{contraction of
$x$ at $p$}. Clearly, $x\arrf{p,I}y$ implies that $x\arr{p,I}y$.
We put $\nu(x)={\displaystyle\sum_{p\in\supp(x)}}\omega^{\hgt(p)}$ (Hessenberg
sum), for any $x\in\nobreak\FFTp$.

We define inductively the relation $\lar{n}$ on $\FFTp$. For $n=0$,
$\lar{n}$ is just the identity relation, while $x\lar{1}y$ iff there
exists $(p,I)\in\MM_T$ such that $x\arr{p,I}y$. Moreover, we put
$x\lar{n+1}y$ whenever there exists $z\in\FFTp$ such that
$x\lar{1}z\lar{n}y$. Furthermore, let $x\lar{*}y$ hold, if
$x\lar{n}y$ for some $n<\omega$. The relations $\larr{n}$ and
$\larr{\ast}$ are defined similarly.

\section{Cancellativity of arrow relations}\label{C:CanArr}

In this section, we fix a colored tree $(T,\utr,\sim)$ and we use the same
notations as in Section~\ref{S:ColTrees}.

\begin{definition}\label{D:CanRel}
A binary relation $\mathbin{R}$ on $\FFTp$ is
\begin{itemize}
\item \emph{additive}, if $x\mathbin{R}y$ implies that
$(x+z)\mathbin{R}(y+z)$, for all $x$, $y$, $z\in\FFTp$.

\item \emph{homogeneous}, if $x\mathbin{R}y$ implies that
$\lambda x\mathbin{R}\lambda y$, for all $x$, $y\in\FFTp$ and
$\lambda\in\FF^+$.

\item \emph{cancellative}, if $(x+z)\mathbin{R}(y+z)$
implies that $x\mathbin{R}y$, for all $x$, $y$, $z\in\FFTp$.
\end{itemize}
\end{definition}

The proof of the following lemma is trivial.

\begin{lemma}\label{L:BasicRelAdd}
The relations $\arr{p,I}$, $\lar{1}$, $\lar{n}$, and $\lar{*}$ are
additive and homogeneous, for all $n<\omega$ and $(p,I)\in\MM_T$.
\end{lemma}

The following lemma states that under certain conditions, arrows of the form
$\arr{p,I}$ may commute.

\begin{lemma}\label{L:CommArr}
Let $x$, $y$, $z\in\FFTp$, let $(p,I)$, $(q,J)\in\MM_T$. If
$x\arr{p,I}y\arr{q,J}z$ and $p\notin J$, then there exists $y'\in\FFTp$
such that $x\arr{q,J}y'\arr{p,I}z$.
\end{lemma}

\begin{proof}
There are $\lambda$, $\mu\in\FF^+$ and $u$, $v\in\FFTp$ such that the
following equalities hold:
       \begin{align}
       x&=\frac{\lambda}{|I|}\sum_{p'\in I}\pd'+u,\label{Eq:zmuq'J}\\
       y&=\lambda\pd+u=\frac{\mu}{|J|}\sum_{q'\in J}\qd'+v,\label{Eq:ylamu}\\
       z&=\mu\qd+v.\label{Eq:xlamu}
       \end{align}
{}From~\eqref{Eq:ylamu} and the assumption that $p\notin J$ it follows that
there exists $w\in\FFTp$ such that
       \[
       u=\frac{\mu}{|J|}\sum_{q'\in J}\qd'+w\text{ and }v=\lambda\pd+w.
       \]
By~\eqref{Eq:zmuq'J} and~\eqref{Eq:xlamu},
$x=\frac{\lambda}{|I|}\sum_{p'\in I}\pd'+\frac{\mu}{|J|}\sum_{q'\in J}\qd'+w$
while $z=\lambda\pd+\mu\qd+w$,
whence $x\arr{q,J}y'\arr{p,I}z$ with
$y'=\frac{\lambda}{|I|}\sum_{p'\in I}\pd'+\mu\qd+w$.
\end{proof}

Now we reach the main result of this section.

\begin{proposition}\label{P:ArrCan}
The relations $\lar{n}$, for $n<\omega$,
and $\lar{*}$ are cancellative.
\end{proposition}

\begin{proof}
It suffices to prove that $\lar{n}$ is cancellative. We argue by induction
on $n$. The statement is trivial for $n=0$. Consider the case where $n=1$.
Since $\lar{1}$ is homogeneous (see Lemma~\ref{L:BasicRelAdd}), it suffices
to prove that $\pd+x\arr{q,I}\pd+y$ implies that $x\arr{q,I}y$, for all
$x$, $y\in\FFTp$, $p\in T$, and $(q,I)\in\MM_T$. By assumption, there are
$\lambda\in\FF^+$ and $u\in\FFTp$ such that
       \begin{align}
       \pd+x&=\frac{\lambda}{|I|}\sum_{r\in I}\rd+u,\label{Eq:pdyru}\\
       \pd+y&=\lambda\qd+u.\label{Eq:pdxqu}
       \end{align}
If $p\neq q$, then, by~\eqref{Eq:pdxqu}, there exists $v\in\FFTp$ such that
$u=\pd+v$. If $p=q$, then $p\notin I$, thus, by~\eqref{Eq:pdyru}, $u=\pd+v$ for
some $v\in\FFTp$. In both cases, $x\arr{q,I}y$. This concludes the $n=1$ case.

Now suppose that $n>1$ and that we have proved the statement for $n-1$.
Let $\pd+x\lar{n}\pd+y$, with $p\in T$ and $x$, $y\in\FFTp$, we
prove that $x\lar{n}y$. There exists $z\in\FFTp$ such that
$\pd+x\lar{n-1}z\lar{1}\pd+y$.

Let $z=z(p)\pd+z'$ where $z'\in\FFTp$ and $z'(p)=0$. Thus, by the
induction hypothesis, either $x\lar{n-1}(z(p)-1)\pd+z'\lar{1}y$ in case
$z(p)\ges 1$, or $(1-z(p))\pd+x\lar{n-1}z'\lar{1}(1-z(p))\pd+y$ in case
$z(p)<1$. In the first case, $x\lar{n}y$, and we are done.
Hence we may assume that
$\pd+x\lar{n-1}z\lar{1}\pd+y$ with $z(p)=0$. {}From $z\lar{1}\pd+y$ and
$z(p)=0$ it follows that $z\arr{p,I}\pd+y$ for some $I\in\MM_T(p)$.
Since $\pd+x\lar{n-1}z$, there exists a chain of the form
       \[
\pd+x=z_0\arr{p_1,I_1}z_1\arr{p_2,I_2}\cdots\arr{p_{n-1},I_{n-1}}z_{n-1}=z,
       \]
where $(p_1,I_1)$, \dots, $(p_{n-1},I_{n-1})\in\MM_T$. Since
$z_0(p)\ges 1>0$ and $z_{n-1}(p)=0$, the largest element $k$ of $\fso{n-1}$
such that $z_k(p)>0$ exists and $k<n-1$. {}From
$z_k\arr{p_{k+1},I_{k+1}}z_{k+1}$ and $z_{k+1}(p)=0$ it follows that
$p_{k+1}=p_*$, in particular,
$\hgt(p_{k+1})<\hgt(p)$. Let $l$ be the largest element of $\fsi{n-1}$ with
$\hgt(p_l)$ minimum; so $\hgt(p_l)<\hgt(p)$. By repeatedly
applying Lemma~\ref{L:CommArr} throughout the chain
       \[
       z_{l-1}\arr{p_l,I_l}z_l\arr{p_{l+1},I_{l+1}}\cdots
       \arr{p_{n-1},I_{n-1}}z_{n-1}=z\arr{p,I}\pd+y,
       \]
(observe that $p_l\notin I_{l+1}\cup\cdots\cup I_{n-1}\cup I$), we obtain a
chain of the form
       \[
z_{l-1}\arr{p_{l+1},I_{l+1}}z'_l\arr{p_{l+2},I_{l+2}}\cdots\arr{p,I}
z'_{n-1}\arr{p_l,I_l}\pd+y,
       \]
with $z'_l$, \dots, $z'_{n-1}\in\FFTp$. Hence, $\pd+x\lar{n-1}z'_{n-1}$.
Furthermore, from $z'_{n-1}\arr{p_l,I_l}\pd+y$ and $\hgt(p_l)<\hgt(p)$ it
follows
that $z'_{n-1}(p)\ges 1$, thus there exists $u\in\FFTp$ such that
$z'_{n-1}=\pd+u$. Hence $\pd+x\lar{n-1}\pd+u\lar{1}\pd+y$, whence, by the
induction hypothesis,
$x\lar{n-1}u\lar{1}y$, thus $x\lar{n}y$.
\end{proof}

\section{Confluence of $\lar{1}$; the relation $\equiv$}
\label{S:Confl}

\begin{lemma}\label{L:RevConf1}
Let $u$, $v$, $x\in\FFTp$. If $x\lar{1}u$ and $x\lar{1}v$, then there exists
$w\in\FFTp$ such that $u\lar{1}w$ and $v\lar{1}w$.
\end{lemma}

\begin{proof}
There are $\lambda$, $\mu\in\FF^+$, $(p,I)$, $(q,J)\in\MM_T$, and $u'$,
$v'\in\FFTp$ such that, putting $m=|I|$ and $n=|J|$, the following
inequalities hold:
       \begin{align}
       u&=\lambda\pd+u',\label{Eq:ulampu'}\\
       v&=\mu\qd+v',\label{Eq:vmuqv'}\\
       x&=\frac{\lambda}{m}\sum_{p'\in I}\pd'+u'=
       \frac{\mu}{n}\sum_{q'\in J}\qd'+v'.\label{Eq:xp'u'q'v'}
       \end{align}
Without loss of generality, $\lambda\les\mu$.
We separate cases.

\Case{1.} $I=J$. Since $T$ is a tree, $p=q$. {}From~\eqref{Eq:xp'u'q'v'}
follows that $u'=\frac{\mu-\lambda}{m}\sum_{p'\in I}\pd'+v'$, thus
       \[
       u=\lambda\pd+\frac{\mu-\lambda}{m}\sum_{p'\in I}\pd'+v'\text{ and }
       v=\mu\pd+v'=\lambda\pd+(\mu-\lambda)\pd+v',
       \]
hence $u\lar{1}v$, so $w=v$ is as desired.

\Case{2.} $I\neq J$. Since both $I$ and $J$ are $\sim$-equivalence classes,
they are disjoint, thus, by~\eqref{Eq:xp'u'q'v'}, there exists $t\in\FFTp$
such that
       \[
       u'=\frac{\mu}{n}\sum_{q'\in J}\qd'+t\text{ and }
       v'=\frac{\lambda}{m}\sum_{p'\in I}\pd'+t,
       \]
whence, by~\eqref{Eq:ulampu'} and~\eqref{Eq:vmuqv'},
       \[
       u=\lambda\pd+\frac{\mu}{n}\sum_{q'\in J}\qd'+t\text{ and }
       v=\frac{\lambda}{m}\sum_{p'\in I}\pd'+\mu\qd+t,
       \]
therefore, $u\lar{1}w$ and $v\lar{1}w$ where $w=\lambda\pd+\mu\qd+t$.
\end{proof}

Now an easy induction proof yields immediately the following lemma.

\begin{lemma}\label{L:RevConfmn}
Let $u$, $v$, $x\in\FFTp$, let $m$, $n<\omega$.
\begin{enumerate}
\item If $x\lar{m}u$ and $x\lar{n}v$, then there exists $w\in\FFTp$ such that
$u\lar{n}w$ and $v\lar{m}w$.

\item If $x\lar{*}u$ and $x\lar{*}v$, then there exists $w\in\FFTp$ such
that $u\lar{*}w$ and $v\lar{*}w$.
\end{enumerate}
\end{lemma}

\begin{notation}\label{No:DefEquiv}
For $x$, $y\in\FFTp$, let $x\equiv y$ hold, if there exists $u\in\FFTp$ such
that $x\lar{*}u$ and $y\lar{*}u$.
\end{notation}

As an immediate consequence of Lemma~\ref{L:RevConfmn}(ii), we obtain the
following.

\begin{proposition}\label{P:EquivTrans}
The relation $\equiv$ is an equivalence relation on $\FFTp$.
\end{proposition}

\section{The relations $\larr{n}$ and the element $\nor{x}$}\label{S:xsharp}

We shall now make use of the relations $\larr{n}$ introduced in
Section~\ref{S:ColTrees}.

\begin{lemma}\label{L:ShnuDecr}
If $x\larr{*}y$ and $x\neq y$, then $\nu(x)>\nu(y)$, for all $x$, $y\in\FFTp$.
\end{lemma}

\begin{proof}
It suffices to consider the case where $x\larr{1}y$. There are decompositions
of the form
       \[
       x=\frac{\lambda}{|I|}\sum_{q\in I}\qd+u\text{ and }y=\lambda\pd+u,
       \]
where $\lambda\in\FF^{++}$, $(p,I)\in\MM_T$, $u\in\FFTp$, and $u(q_0)=0$
for some $q_0\in I$. It follows that
       \[
       \supp(x)=\supp(u)\cup I\text{ while }\supp(y)=\supp(u)\cup\set{p},
       \]
with $q_0\in I\setminus\supp(u)$. Therefore, using again the Hessenberg
addition,
       \begin{align}
       \nu(y)&\leq\nu(u)+\omega^{\hgt(p)}\notag\\
       &<\nu(u)+\omega^{\hgt(q_0)}&&(\text{because }\hgt(q_0)=\hgt(p)+1)
       \notag\\
       &\leq\nu(x)&&(\text{because }u(q_0)=0).\tag*{\qed}
       \end{align}
\renewcommand{\qed}{}
\end{proof}

\begin{lemma}\label{L:expIpred}
For all $x\in\FFTp$ and all $(p,I)\in\MM_T$, there exists $y\in\FFTp$
such that $x\arrf{p,I}y$.
\end{lemma}

\begin{proof}
We have $x=\frac{\lambda}{|I|}\sum_{q\in I}\qd+u$,
where $\lambda=|I|\cdot\min\setm{x(q)}{q\in I}$ and
$u=x-\frac{\lambda}{|I|}\sum_{q\in I}\qd$.
Take $y=\lambda\pd+u$. Obviously, $x\arrf{p,I}y$.
\end{proof}

\begin{lemma}\label{L:TriExpan}
Let $x$, $y$, $z\in\FFTp$, let $(p,I)\in\MM_T$. If $z\arrf{p,I}x$ and
$z\arr{p,I}y$, then $y\arrf{p,I}x$. Furthermore, $y\neq z$ implies that
$x\neq z$.
\end{lemma}

\begin{proof}
There are $\lambda$, $\mu\in\FF^+$ and $u$, $v\in\FFTp$ such that, putting
$n=|I|$, the following equalities hold:
       \begin{align}
       x&=\lambda\pd+u,\label{Eq:xlampu}\\
       y&=\mu\pd+v,\label{Eq:ylampv}\\
       z&=\frac{\lambda}{n}\sum_{q\in I}\qd+u=\frac{\mu}{n}\sum_{q\in I}\qd+v
       \label{Eq:zlamuuv}
       \end{align}
with $u(q_0)=0$ for some $q_0\in I$. Thus, by~\eqref{Eq:zlamuuv},
$\mu\les\lambda$, whence $v=\frac{\lambda-\mu}{n}\sum_{q\in I}\qd+u$.
Therefore,
       \[
       x=\lambda\pd+u\text{ and
}y=\mu\pd+\frac{\lambda-\mu}{n}\sum_{q\in I}\qd+u,
       \]
with $u(q_0)=0$, whence $y\arrf{p,I}x$.

If $y\neq z$, then $\mu>0$, thus $\lambda>0$, thus $x\neq z$.
\end{proof}

\begin{definition}\label{D:origin1}
For $x\in\FFTp$, we put
       \begin{align*}
       \Phi(x)&=\setm{y\in\FFTp}{x\lar{*}y},\\
       \Phi^*(x)&=\setm{y\in\Phi(x)}{\Phi(y)=\set{y}}.
       \end{align*}
\end{definition}

\begin{lemma}\label{L:Phi*xsingl}
The set $\Phi^*(x)$ is a singleton, for all $x\in\FFTp$.
\end{lemma}

\begin{proof}
Let $u$ be an element of $\Phi(x)$ with $\nu(u)$ smallest possible.
Suppose that there exists $v\neq u$ such that $u\lar{*}v$. Then there
exists $v\neq u$ such that $u\lar{1}v$, thus, by Lemmas~\ref{L:expIpred}
and \ref{L:TriExpan}, there exists $v\neq u$ such that $u\larr{1}v$. By
Lemma~\ref{L:ShnuDecr}, $\nu(v)<\nu(u)$, which contradicts the minimality
assumption on $\nu(u)$. Therefore, $u$ belongs to $\Phi^*(x)$. The
uniqueness statement on $u$ follows from Lemma~\ref{L:RevConfmn}.
\end{proof}

\begin{definition}\label{D:origin}
Let the \emph{normal form} of $x\in\FFTp$ be the unique element of
$\Phi^*(x)$; we denote it by $\nor{x}$. We say that $x$ is
\emph{normal}, if $x=\nor{x}$.
\end{definition}

Therefore, $x\lar{*}y$ implies that $\nor{x}=\nor{y}$.

We leave to the reader the easy proof of the following lemma.

\begin{lemma}\label{L:somenormal}\hfill
\begin{enumerate}
\item Every element of the form $\lambda\pd$, where $\lambda\in\FF^+$ and
$p\in T$, is normal.

\item If $x$ is a normal element, then any $y\in\FFTp$ such that $y\leq x$ is
normal.
\end{enumerate}
\end{lemma}

\begin{remark}\label{Rk:Antis}
It can be proved that the relation $\lar{*}$ is \emph{antisymmetric}.
However, we will not use this fact.
\end{remark}

Now we are coming to the main result of this section.

\begin{lemma}\label{L:ShDa}
If $x\lar{n}\nor{x}$, then $x\larr{n}\nor{x}$, for all $x\in\FFTp$ and
$n<\omega$.
\end{lemma}

\begin{proof}
We argue by induction on $n$. If $n=0$ then $x=\nor{x}$, and we are done.
Suppose that $n>0$. There exist $y\in\FFTp$
and $(p,I)\in\MM_T$ such that $x\arr{p,I}y\lar{n-1}\nor{x}$. By
Lemma~\ref{L:expIpred}, there exists $z\in\FFTp$ such that $x\arrf{p,I}z$.
Now, by Lemma~\ref{L:TriExpan}, $y\arrf{p,I}z$.
By Lemma~\ref{L:RevConfmn}, there exists $w\in\FFTp$ such that
$\nor{x}\lar{1}w$ and $z\lar{n-1}w$.
Thus, $w=\nor{x}$ and $z\lar{n-1}\nor{x}$. Since $\nor{x}=\nor{z}$, we get,
by the induction hypothesis, that $x\larr{1}z\larr{n-1}\nor{x}$.
\end{proof}

\section{The cancellation theorem}\label{S:CanThm}
We first establish a technical lemma.

\begin{lemma}\label{L:larrtolamp+x}
Let $x$, $y\in\FFTp$ with $x$ normal, let $p\in T$, let $\lambda\in\FF^+$.
If $\lambda\pd+x\larr{1}y$ and $\lambda\pd+x\neq y$, then $\lambda>0$,
$x(p)=0$, and there are $x'\in\FFTp$ and $\xi\in(0,\lambda]$ in $\FF$ such
that, putting $I=[p]\setminus\set{p}$ and $l=|I|$, the following
statements hold:
\begin{enumerate}
\item $x=\xi\sum_{q\in I}\qd+x'$ and
$y=(\lambda-\xi)\pd+(l+1)\xi\pd_*+x'$.

\item $(\lambda-\xi)\pd+x'$ is normal.

\item $((\lambda-\xi)\pd+x')(q)=x(q)$, for all $q\in T$ such that
$\hgt(p)<\hgt(q)$.
\end{enumerate}
\end{lemma}

\begin{proof}
If $\lambda=0$, then, since $x$ is normal, $y=x$, \contr; whence $\lambda>0$.
Suppose now that $x(p)>0$. Then there exists $\eps\in\FF^{++}$
such that $\eps\lambda\pd\leq(1-\eps)x$, thus $\eps(\lambda\pd+x)\leq x$. Since
$x$ is normal, by Lemma~\ref{L:somenormal}(ii), $\lambda\pd+x$ is normal,
\contr\ with $\lambda\pd+x\larr{1}y$ and $\lambda\pd+x\neq y$.
Hence $x(p)=0$.

Put $I=\setm{p_i}{0\les i<l}$.
Let $\xi'$ be the least element of
$\setm{x(p_i)}{i<l}$, and put $\xi=\min\set{\xi',\lambda}$. Since $x$ is
normal, the contraction from $\lambda\pd+x$ to $y$ occurs at $p_*$, and
there are decompositions of the form
       \begin{align}
       \lambda\pd+x&=(\lambda-\xi)\pd+\xi\pd+\xi\sum_{i<l}\pd_i+x',
       \label{Eq:decpdx}\\
       y&=(\lambda-\xi)\pd+(l+1)\xi\pd_*+x',\label{Eq:decy}
       \end{align}
with $x'\in\FFTp$, $x'(p_j)=0$ for some $j<l$, and, since $y\neq\lambda\pd+x$,
$\xi>0$.

The element $(\lambda-\xi)\pd+x'$ is normal, otherwise, $\xi<\lambda$, and,
by the same argument as in the previous paragraph, there exists
$\xi''\in\FF^{++}$ such that $x'\geq\xi''\sum_{i<l}\pd_i$, which contradicts
$x'(p_j)=0$.

For $q\in T$ such that $\hgt(p)<\hgt(q)$, it follows from~\eqref{Eq:decpdx}
that $x(q)=x'(q)$, whence $((\lambda-\xi)\pd+x')(q)=x(q)$.
\end{proof}

Now we are ready to prove the main result of this section.

\begin{theorem}[The cancellation theorem]\label{T:EquivCan}
The relation $\equiv$ is cancellative.
\end{theorem}

\begin{proof}
We recall that $\equiv$ is an equivalence relation (see
Proposition~\ref{P:EquivTrans}). Observe that by Lemma~\ref{L:BasicRelAdd},
$\equiv$ is additive.

We need to prove that for all $x$, $y\in\FFT$, and all $p\in T$, if
$\pd+x\equiv\pd+y$, then $x\equiv y$. Since $x\equiv\nor{x}$ and
$y\equiv\nor{y}$, it suffices to consider the case where both $x$ and $y$ are
normal, and then $\pd+x\lar{*}u$ and $\pd+y\lar{*}u$ where
$u=\nor{(\pd+x)}=\nor{(\pd+y)}$. Thus it suffices to prove the following
statement:
       \begin{align*}
       &\text{For all }m,\,n<\omega,\ p\in T,\text{ and }x,\,y,\,u\in\FFTp
       \text{ normal,}\\
       &\text{if }\pd+x\lar{m}u\text{ and }\pd+y\lar{n}u,\text{ then }x=y.
       \end{align*}
We argue by induction on $m+n$. If $m=0$, then $\pd+y\lar{n}\pd+x$, thus, by
Proposition~\ref{P:ArrCan}, $y\lar{n}x$, thus, since $y$ is normal, $x=y$, so
we are done. A similar argument holds if $n=0$.

Suppose from now on that $m$ and $n$ are nonzero. It follows from
Lemma~\ref{L:ShDa} that $\pd+x\larr{m}u$ and $\pd+y\larr{n}u$, thus there are
chains of the form
       \begin{align}
       \pd+x&=x_0\larr{1}x_1\larr{1}\cdots\larr{1}x_{m-1}\larr{1}x_m=u,
       \label{Eq:chu2px}\\
       \pd+y&=y_0\larr{1}y_1\larr{1}\cdots\larr{1}y_{n-1}\larr{1}y_n=u,
       \label{Eq:chu2py}
       \end{align}
for some $x_0$,\ldots, $x_m$, $y_0$,\ldots, $y_n\in\FFTp$.
If two distinct occurrences of one of the chains~\eqref{Eq:chu2px},
\eqref{Eq:chu2py} are equal, then either $\pd+x\larr{m-1}u$ or
$\pd+y\larr{n-1}u$, thus $x=y$ by the induction hypothesis.

Suppose from now on that each of the chains~\eqref{Eq:chu2px} and
\eqref{Eq:chu2py} has all its entries distinct. Let $(q_i)_{i<l}$ be a
one-to-one enumeration of $[p]\setminus\set{p}$. By Lemma~\ref{L:larrtolamp+x},
there are $\xi$, $\eta\in(0,1]$ in $\FF$, together with
$\ol{x}$, $\ol{y}\in\FFTp$, $i_0<m$, and $j_0<n$ such that the following
equations hold:
       \begin{align}
  \pd+x&=(1-\xi)\pd+\xi\pd+\xi\sum_{i<l}\qd_i+\ol{x},\label{Eq:DecLastLev1}\\
       x_1&=(1-\xi)\pd+(l+1)\xi\pd_*+\ol{x},\label{Eq:DecLastLev2}\\
     \pd+y&=(1-\eta)\pd+\eta\pd+\eta\sum_{i<l}\qd_i+\ol{y},
     \label{Eq:DecLastLev3}\\
       y_1&=(1-\eta)\pd+(l+1)\eta\pd_*+\ol{y},\label{Eq:DecLastLev4}\\
       \ol{x}(q_{i_0})&=\ol{y}(q_{j_0})=0.\label{Eq:MaxDec}
       \end{align}
Furthermore, by Lemma~\ref{L:larrtolamp+x}, $x(p)=y(p)=0$ and both elements
$x'_1=(1-\xi)\pd+\ol{x}$ and $y'_1=(1-\eta)\pd+\ol{y}$ are normal.
Observe that
       \begin{equation}\label{Eq:xyfromx'y'}
       x_1=(l+1)\xi\pd_*+x'_1\text{ and }y_1=(l+1)\eta\pd_*+y'_1.
       \end{equation}
Define inductively $p_0=p$, and $p_{i+1}=(p_i)_*$ (for $i<\omega$) whenever it
is defined. In particular, $p_1=p_*$. By applying inductively
Lemma~\ref{L:larrtolamp+x}, starting with~\eqref{Eq:DecLastLev2}, we obtain
decompositions $x_i=\lambda_i\pd_i+x'_i$, for $1\les i\les m$, with
$\lambda_i\in\FF^{++}$ and $x'_i\in\FFTp$ normal such that
$x'_i(p)=x'_1(p)=1-\xi$. Similarly, starting with~\eqref{Eq:DecLastLev4}, we
obtain decompositions $y_j=\mu_j\pd_j+y'_j$, for
$1\les j\les n$, with $\mu_j\in\FF^{++}$ and $y'_j\in\FFTp$ normal such that
$y'_j(p)=y'_1(p)=1-\eta$.

In particular, $1-\xi=x_m(p)=u(p)=y_n(p)=1-\eta$, whence $\xi=\eta$. Hence,
$x_1=(l+1)\xi\pd_*+x'_1$ and $y_1=(l+1)\xi\pd_*+y'_1$ with
both $x'_1$ and $y'_1$ normal, thus, since $x_1\lar{m-1}u$ and
$y_1\lar{n-1}u$ and by the induction hypothesis, $x'_1=y'_1$; whence
$\ol{x}=\ol{y}$. Therefore, by~\eqref{Eq:DecLastLev1} and
\eqref{Eq:DecLastLev3},
$x=y$, which concludes the proof.
\end{proof}

\begin{notation}\label{No:N(T)}
Let $N_T$ denote the subspace of $\FFT$ defined by
       \[
       N_T=\setm{x-y}{x,\,y\in\FFTp\text{ and }x\equiv y}.
       \]
We put $\VV_T=\FFT/{N_T}$, and we put $\ol{p}=\pd+N_T$, for all $p\in T$.
\end{notation}

As an immediate consequence of Theorem~\ref{T:EquivCan}, we obtain the
following.

\begin{corollary}\label{C:EquivCan}
For all $x$, $y\in\FFTp$, $x-y\in N_T$ if{f} $x\equiv y$.
\end{corollary}

In particular, since all elements $\pd$, for $p\in T$, are normal, we obtain:

\begin{corollary}\label{C:olpdist}
The map $p\mapsto\ol{p}$ is one-to-one.
\end{corollary}

\section{Plenary subsets, plenary embeddings, and the trace functional}
\label{S:Plen}

\begin{definition}\label{D:PlSub}
A subset $\Omega$ of a vector space $V$ is \emph{plenary}, if for every
$x\in\Col(\Omega)$, there exists a least
(necessarily finite) subset $X\in\Co(V,\Omega)$ such
that $x\in\Col(X)$.
\end{definition}

Observe that for a subset
$\Omega$ of $V$, the canonical map
$\varphi_\Omega\colon\Co(V,\Omega)\to\Co(V)$, $X\mapsto\Col(X)$ is always a
complete join-embedding. We leave to the reader the straightforward
proof of the following.

\begin{proposition}\label{P:Plenary}
Let $\Omega$ be a subset of a vector space $V$ over a totally ordered
division ring. Then $\Omega$ is plenary if{f} the canonical
map $\varphi_\Omega$ from $\Co(V,\Omega)$ into $\Co(V)$ is a complete lattice
embedding.
\end{proposition}

\begin{example}\label{Ex:PlnonPl}
The whole space $V$, or any affinely independent subset of $V$, is plenary.
On the other hand, the square $C=\set{(0,0),(0,1),(1,0),(1,1)}$ is
not plenary in
$\QQ^2$ (take $X=C\setminus\set{(1,1)}$, $Y=C\setminus\set{(1,0)}$).
\end{example}

\begin{definition}\label{D:PlEmb}
For a join-semilattice $L$ and a vector space $V$ over a totally ordered
division ring, a map $\varphi\colon L\to\Co(V)$ is \emph{plenary}, if
$\varphi=\varphi_\Omega\circ\psi$ for some plenary subset~$\Omega$ of~$V$
and some join-homomorphism $\psi\colon L\to\Co(V,\Omega)$ that preserves
existing meets.
\end{definition}

Hence every plenary map from a lattice to $\Co(V)$ is a lattice
homomorphism, and it preserves existing meets. Furthermore, in the statement
above, $\varphi$ is an embedding if{f} $\psi$ is an embedding.

{}From now on until the end of the present section, we shall fix a totally
ordered division ring $\FF$ and a colored tree
$(T,\utr,\sim)$. We shall use the notations and terminology of the previous
sections about $\FFT$, $\lar{*}$, $\equiv$, $\VV_T$, $N_T$, $\pd$,
$\ol{p}$, and
so on.

\begin{lemma}\label{L:trace}
There exists a unique linear functional $\tau\colon\VV_T\to\FF$ such that
$\tau(\ol{p})=1$ for all $p\in T$.
\end{lemma}

\begin{proof}
Let $f\colon\FFT\to\FF$ be the unique linear functional defined by
$f(\pd)=1$ for
all $p\in T$. It is sufficient to prove that the restriction of $f$ to $N_T$ is
zero. For this, it is sufficient to prove that $x\lar{1}y$ implies that
$f(x)=f(y)$, for all $x$, $y\in\FFTp$, which is obvious.
\end{proof}

We shall call the \emph{trace functional} the linear functional
$\tau\colon\VV_T\to\FF$ given by Lemma~\ref{L:trace}.

\begin{notation}\label{No:OT}
Set $\Omega_T=\setm{\ol{p}}{p\in T}$, a subset of $\VV_T$.
For $x\in\FFTp$, we set $\bsupp(x)=\setm{\ol{p}}{p\in\supp(x)}$, a subset of
$\Omega_T$.
\end{notation}

\begin{lemma}\label{L:ConvPres}
Let $x$, $y\in\FFTp$. If $x\lar{*}y$, then
$\bsupp(y)\subseteq\Omega_T\cap\Col(\bsupp(x))$.
\end{lemma}

\begin{proof}
It suffices to verify this for $x\lar{1}y$ and $x\neq y$. There are
$\lambda\in\FF^{++}$, $(p,I)\in\MM_T$, and $u\in\FFTp$ such that
$x=\frac{\lambda}{|I|}\sum_{q\in I}\qd+u$ and $y=\lambda\pd+u$, whence
$\bsupp(y) =\bsupp(u)\cup\set{\ol{p}}$ while
$\bsupp(x)=\bsupp(u)\cup\setm{\ol{q}}{q\in I}$. Hence
$\ol{p}=\frac{1}{|I|}\sum_{q\in I}\ol{q}$ belongs to $\Col(\bsupp(x))$.
\end{proof}

\begin{proposition}\label{P:OmPlen}
The set $\Omega_T$ is a plenary subset of $\VV_T$.
\end{proposition}

\begin{proof}
Let $\boldsymbol{x}\in\Col(\Omega_T)$. Denote
$\ol{Y}=\setm{\ol{p}}{p\in Y}$, for all $Y\subseteq T$, and denote by $x$
the unique normal representative of $\boldsymbol{x}$.
There are a positive integer $m$, scalars $\alpha_0$, \dots,
$\alpha_{m-1}\in\FF^{++}$, and elements $p_0$, \dots, $p_{m-1}\in T$ such
that
      \begin{equation}\label{Eq:NormFmx}
      x=\sum_{i<m}\alpha_i\pd_i.
      \end{equation}
{}From $\boldsymbol{x}\in\Col(\Omega_T)$ and Lemma~\ref{L:trace} it follows
that $\tau(\boldsymbol{x})=1$, that is, $\sum_{i<m}\alpha_i=1$. Hence, by
\eqref{Eq:NormFmx}, $x\in\Col(\ol{X})$, where we put $X=\setm{p_i}{i<n}$.

Let $Y\subseteq T$ such that $\boldsymbol{x}\in\Col(\ol{Y})$. There are a
positive integer $n$, scalars $\beta_0$, \dots, $\beta_{n-1}\in\FF^{++}$,
and elements $q_0$, \dots, $q_{m-1}\in Y$ such that
      \begin{equation}\label{Eq:DecFmx}
      \boldsymbol{x}=\sum_{j<n}\beta_j\ol{q}_j.
      \end{equation}
Put $y=\sum_{j<n}\beta_j\qd_j$. It follows from~\eqref{Eq:NormFmx} and
\eqref{Eq:DecFmx} that $x\equiv y$, but $x$ is normal, thus $y\lar{*}x$.
By Lemma~\ref{L:ConvPres},
      \[
      \ol{X}=\bsupp(x)\subseteq\Omega_T\cap\Col(\bsupp(y))\subseteq\ol{Y},
      \]
which proves that $\ol{X}$ is the least subset of $\Omega_T$ whose convex
hull contains $\boldsymbol{x}$.
\end{proof}

By Proposition~\ref{P:JirrCoVO}, the \jirr\ elements of $\Co(\VV_,\Omega_T)$
are the trivial ones. We obtain another remarkable property of the set
$\Omega_T$.

\begin{proposition}\label{P:noDcycle}
For all $p$, $q\in T$, $\set{\ol{p}}\DD\set{\ol{q}}$ in
$\Co(\VV_T,\Omega_T)$ implies that $p\tr q$. In particular, the
join-dependency relation of $\Co(\VV_T,\Omega_T)$ is well-founded
\pup{i.e., it has no infinite descending sequence} on the set of
\jirr\ elements of $\Co(\VV_T,\Omega_T)$.
\end{proposition}

\begin{proof}
Since $p\tr q$
implies that $\hgt(p)<\hgt(q)$, it suffices to prove the first assertion. By
assumption, $p\neq q$ and there exists $X\in\Co(\VV_T,\Omega_T)$ such that
$\ol{p}\notin X$ and $\ol{p}\in\set{\ol{q}}\vee X$, thus there are
$\lambda\in\FF$ with $0<\lambda<1$ and $\boldsymbol{x}\in\Col(X)$ such that
       \[
       \ol{p}=(1-\lambda)\ol{q}+\lambda\boldsymbol{x}.
       \]
Since $\pd$ is normal and by Corollary~\ref{C:EquivCan}, it follows that
       \[
       (1-\lambda)\qd+\lambda x\lar{*}\pd,
       \]
for some (any) $x\in\boldsymbol{x}$.
In particular, from $p\neq q$ it follows that $p\tr q$.
\end{proof}

As, in the finite case, the nonexistence of $\DD$-cycles is equivalent to being
lower bounded (see~\cite{FJN}), we obtain the following.

\begin{corollary}\label{C:noDcycle}
If $T$ is finite, then $\Co(\VV_T,\Omega_T)$ is finite lower bounded.
\end{corollary}

\section{Norms on trees}\label{S:Norms}

\begin{definition}\label{D:Norm}
Let $T$ be a colored tree, let $L$ be a join-semilattice. A \emph{$L$-valued
norm} on~$T$ is a map $e\colon T\to L^-$ which satisfies the following
conditions:
\begin{enumerate}
\item For all $(p,I)\in\MM_T$, $e[I]=\setm{e(q)}{q\in I}$ is a nontrivial
join-cover of $e(p)$.

\item For all $p\in T$ and every nontrivial join-cover $X$ of $e(p)$, there
exists $I\in\MM_T(p)$ such that $e[I]\ll X$.
\end{enumerate}

In addition, we say that $e$ is \emph{full}, if every element $x$ of $L$ is the
join of all elements of~$e[T]$ below $x$.
\end{definition}

The main goal of this section is to prove the following result.

\begin{theorem}\label{T:Norm2Emb}
Let $T$ be a colored tree, let $L$ be a join-semilattice, let
$e\colon T\to L^-$ be a norm, let $\FF$ be a totally ordered
division ring. Consider the vector space $\VV_T$ and the subset $\Omega_T$
constructed in previous sections from $T$ and $\FF$. Then one can define a
join-homomorphism $\psi\colon L\to\Co(\VV_T,\Omega_T)$ by the rule
       \[
       \psi(x)=\setm{\ol{p}}{p\in T\text{ and }e(p)\leq x},\text{ for
all }x\in L.
       \]
Then
$\psi$ preserves existing meets.
Furthermore, the following statements hold:
\begin{enumerate}
\item The map $\varphi\colon L\to\Co(\VV_T)$ defined by
$\varphi(x)=\Col(\psi(x))$, for all $x\in L$, is a plenary join-homomorphism
from $L$ to $\Co(\VV_T)$.

\item Both $\psi$ and $\varphi$ are zero-preserving.

\item If the norm $e$ is full, then both $\psi$ and $\varphi$ are
embeddings.
\end{enumerate}
\end{theorem}

\begin{proof}
Put $L^\circ=L\cup\set{\mathrm{O}}$, for a new zero element $\mathrm{O}$.
We first extend $e$ to a map from $\FFTp$ to $L^\circ$, still denoted by $e$,
as follows:
       \[
       e(x)=\bigvee\setm{e(p)}{p\in\supp(x)},\text{ for all }x\in\FFTp,
       \]
with the convention $\bigvee\es=\mathrm{O}$.

\setcounter{claim}{0}
\begin{claim}\label{Cl:edecr}
If $x\lar{*}y$, then $e(y)\leq e(x)$, for all $x$, $y\in\FFTp$.
\end{claim}

\begin{scproof}
It suffices to prove the result in the case where $x\lar{1}y$ and $x\neq y$.
There are $\lambda\in\FF^{++}$, $(p,I)\in\MM_T$, and $z\in\FFTp$ such that
       \[
       x=\frac{\lambda}{|I|}\sum_{q\in I}\qd+z\text{ and }y=\lambda\pd+z.
       \]
Since $e$ is a norm, $e(p)\leq\bigvee e[I]$, whence
       \begin{equation}
       e(y)=e(p)\vee e(z)\leq\bigvee e[I]\vee e(z)=e(x).\tag*{\qed}
       \end{equation}
\renewcommand{\qed}{}
\end{scproof}

\begin{claim}\label{Cl:psixconv}
The set $\psi(x)$ belongs to $\Co(\VV_T,\Omega_T)$, for all $x\in L$.
\end{claim}

\begin{scproof}
Let $p\in T$, suppose that $\ol{p}\in\Col(\psi(x))$, we prove that
$\ol{p}\in\psi(x)$. By assumption, $\ol{p}=\sum_{i<n}\lambda_i\ol{p}_i$ for
some $n>0$, $(\lambda_i)_{i<n}\in\del{n}{\FF}$, $(p_i)_{i<n}\in T^n$ with
$e(p_i)\leq x$, for all $i<n$. By Corollary~\ref{C:EquivCan},
$\pd\equiv\sum_{i<n}\lambda_i\pd_i$, but $\pd$ is normal (see
Lemma~\ref{L:somenormal}(i)), thus $\sum_{i<n}\lambda_i\pd_i\lar{*}\pd$, thus,
by Claim~\ref{Cl:edecr}, $e(p)\leq\bigvee_{i<n}e(p_i)\leq x$, that is,
$\ol{p}\in\psi(x)$.
\end{scproof}

Since $e[T]$ is contained in $L^-$, $\psi(0)=\varphi(0)=\es$ if $L$ has a zero.

It is obvious that $\psi$ preserves existing meets. Now we prove that
$\psi$ is a join-homomorphism.
It is sufficient to prove that for all $x$, $y\in L$ and all $p\in T$, if
$e(p)\leq x\vee y$, then $\ol{p}\in\psi(x)\vee\psi(y)$ (the join
$\psi(x)\vee\psi(y)$ is computed in $\Co(\VV_T,\Omega_T)$). This is obvious if
either $e(p)\leq x$ or $e(p)\leq y$, in which case
$\ol{p}\in\psi(x)\cup\psi(y)$. Suppose that $e(p)\nleq x,y$. Then $\set{x,y}$
is a nontrivial join-cover of $e(p)$, thus, since $e$ is a norm, there exists
$I\in\MM_T(p)$ such that $e[I]\ll\set{x,y}$. Therefore,
$\ol{p}=\frac{1}{|I|}\sum_{q\in I}\ol{q}$ belongs to
$\Col(\psi(x)\cup\psi(y))$,
but $\ol{p}\in\Omega_T$, whence $\ol{p}\in\psi(x)\vee\psi(y)$.

Since $\Omega_T$ is a plenary subset of $\VV_T$ (see
Proposition~\ref{P:OmPlen}), $\varphi$ is a plenary homomorphism.

Finally, suppose that $e$ is a full norm, we prove that $\psi$ is an embedding
(thus $\varphi$ is also an embedding). Let $x$, $y\in L$ such that $x\nleq y$.
Since $e$ is full, there exists $p\in T$ such that $e(p)\leq x$ and
$e(p)\nleq y$, whence $\ol{p}\in\psi(x)\setminus\psi(y)$; thus
$\psi(x)\not\subseteq\psi(y)$. Hence $\varphi$ is an embedding from $L$ into
$\Co(\VV_T)$.
\end{proof}

The result of Theorem~\ref{T:Norm2Emb} for $\FF=\QQ$ does not trivially imply
the result for other totally ordered division rings, as, for example, the
canonical embedding from $\Co(\QQ)$ into $\Co(\RR)$ does not preserve existing
meets.

Although the results of Sections~\ref{S:EmbThms} and \ref{S:RelAbs} are
formulated for lattices, we shall need in subsequent work the semilattice
formulation of Section~\ref{S:Norms}.

\section{Embedding lattices into lattices of convex sets}\label{S:EmbThms}

In this section, we shall apply the results of the previous sections, in order
to represent lattices as lattices of convex sets in vector spaces. Throughout
this section, we shall fix a totally ordered division ring $\FF$.

\begin{theorem}\label{T:ArbLatt}
Every lattice has a plenary, zero-preserving embedding into $\Co(V)$, for some
$\FF$-vector space $V$.
\end{theorem}

\begin{proof}
Let $L$ be a lattice and let $T$ denote the set of all finite sequences
of the form
       \begin{equation}\label{Eq:Exprp}
       p=\seq{a_0,I_0,a_1,I_1,\dots,a_{m-1},I_{m-1},a_m},
       \end{equation}
where $m<\omega$, $a_0$, \dots, $a_m\in L^-$, $I_k$ is a nontrivial join-cover
of $a_k$ and $a_{k+1}\in I_k$, for all $k<m$.
For $p$ given by~\eqref{Eq:Exprp} and $q$ given by
       \begin{equation}\label{Eq:Exprq}
       q=\seq{b_0,J_0,b_1,J_1,\dots,b_{n-1},J_{n-1},b_n},
       \end{equation}
let $p\utr q$ hold, if $p$ is an initial segment of $q$, and let $p\sim q$
hold, if $m=n$ and $(a_k,I_k)=(b_k,J_k)$ for all $k<m$.
Also, let $e(p)=a_m$ if $p$ is given by~\eqref{Eq:Exprp}.
The verification that $(T,\utr,\sim)$ is a colored tree
is straightforward. For $p$ as in~\eqref{Eq:Exprp} and $q$ as
in~\eqref{Eq:Exprq}, $p\prec q$ if{f} $p\utr q$ and $n=m+1$, and then $I=[q]$
consists exactly of those elements of $T$ of the form
       \[
       q'=\seq{a_0,I_0,a_1,I_1,\dots,a_{m-1},I_{m-1},a_m,J_m,x},
       \text{ where }x\in J_m.
       \]
In particular, $e[I]=J_m$ is a nontrivial join-cover of $e(p)=a_m$. As
every nontrivial join-cover of $a_m$ arises in this fashion, $e$ is a
full norm.
\end{proof}

Now for the finite lower bounded case, we get a more precise result. For a
vector space $V$ over a totally ordered division ring, we denote by
$\mathbf{K}(V)$ the lattice of all \emph{convex polytopes} of $V$, that is, the
finitely generated convex subsets of $V$. It is well-known that $\mathbf{K}(V)$
is a \jsd\ sublattice of $\Co(V)$, see Theorem~15 in G. Birkhoff and M.\,K.
Bennett~\cite{BB}.

\begin{theorem}\label{T:FinLB}
Every finite lower bounded lattice $L$ has a plenary, zero-preserving
embedding into $\mathbf{K}(\FF^n)$, for some $n<\omega$.
Furthermore, $L$ has a
zero-preserving embedding
into a lower bounded lattice of the form $\Co(\QQ^n,\Omega)$, for some
$n<\omega$ and some plenary finite subset $\Omega$ of~$\ZZ^n$.
\end{theorem}

\begin{proof}
Let $L$ be a finite lower bounded lattice and let $T$ be the set of all finite
sequences of the form given in
\eqref{Eq:Exprp}, where $n<\omega$, $a_0$, \dots, $a_n\in\J(L)$, $I_k$ is a
minimal nontrivial join-cover of $a_k$ and $a_{k+1}\in I_k$, for all
$k<n$. We define the relations $\utr$ and $\sim$ and the map $e$ as in
the proof of Theorem~\ref{T:ArbLatt}.
The verifications that $(T,\utr,\sim)$ is a colored tree and that
$e$ is a full norm are mostly as in the proof of Theorem~\ref{T:ArbLatt}.
Moreover, since
$L$ is finite lower bounded, it has no $\DD$-cycle, thus $T$ is finite; whence
$\VV_T$ is finite-dimensional and $\Omega_T$ is finite.
By Proposition~\ref{P:OmPlen}, $\Omega_T$ is plenary. By
Corollary~\ref{C:noDcycle}, $\Co(\VV_T,\Omega_T)$ is finite lower bounded. In
case $\FF=\QQ$, fixing an isomorphism from $\VV_T$ onto some
$\QQ^n$ and replacing $\Omega_T$ by $\Omega=m\Omega_T$, for a suitable positive
integer $m$, turns $\Omega$ to a subset of~$\ZZ^n$.

The conclusion follows again from Theorem~\ref{T:Norm2Emb}.
\end{proof}

Hence we have obtained a new universal class of finite lower bounded
lattices, namely, the class of lattices of the form
$\Co(\QQ^n,\Omega)$, where~$n$ is a positive integer and~$\Omega$ is a finite
plenary subset of $\ZZ^n$.
We recall that two other well-known
universal classes of finite lower bounded
lattices consist of the lattices of the form $\mathbf{Sub}_{\wedge}(\two^m)$
(the lattice of all meet-subsemilattices of the Boolean lattice $\two^m$) and
of the lattices of the form $\mathcal{O}(\mathbf{n})$ (the lattice
of all suborders of a given linear order on the finite set~$\mathbf{n}$),
respectively.

\section{The lattices $\Co(V,\Omega)$ and $\mathbf{K}(V)$}\label{S:RelAbs}

In~\cite{AGT}, the problem of embeddability of a given finite lattice into some
finite lattice of the form $\Co(\RR^n,\Omega)$, for a positive integer $n$ and
a finite subset $\Omega$ of $\RR^n$, is posed. The following easy result
establishes a simple relation between embeddability into some $\Co(\FF^n)$
and embeddability into some $\Co(\FF^n,\Omega)$.

\begin{proposition}\label{P:Abs2Rel}
Let $L$ be a lattice, let $V$ be a vector space over a totally ordered
division ring $\FF$, let $\varphi\colon L\into\mathbf{K}(V)$ be a lattice
embedding. Let $\Omega$ be any subset of~$V$ containing all the
extreme points of
all elements of the form $\varphi(x)$, for $x\in L$. Then the map
$\psi\colon L\into\Co(V,\Omega)$, $x\mapsto\varphi(x)\cap\Omega$ is a lattice
embedding from $L$ into $\Co(V,\Omega)$, and $\varphi(x)=\Col(\psi(x))$ for all
$x\in L$.
\end{proposition}

\begin{proof}
It is obvious that $\psi$ is a meet homomorphism. Since every element of the
range of $\varphi$ is the convex hull of its (finite) set of extreme points,
which is contained in~$\Omega$, the equality $\varphi(x)=\Col(\psi(x))$
holds for all $x\in L$, thus $\psi$ is an order-embedding.

Denote by $\eb(X)$ the set of extreme points of a convex polytope $X$ of $V$.
Let $x$, $y\in L$. For any $X\in\Co(V,\Omega)$, if $\psi(x)\vee\psi(y)$ is
contained in $X$, then $\eb(\varphi(x))\cup\eb(\varphi(y))$ is contained in
$X$, thus also the smaller set $\eb(\varphi(x)\vee\varphi(y))$, which is equal
to $\eb(\varphi(x\vee y))$. Hence $\psi(x\vee y)$ is contained in $X$, which
proves that $\psi(x\vee y)=\psi(x)\vee\psi(y)$. Hence $\psi$ is a lattice
homomorphism.
\end{proof}

\begin{corollary}\label{C:Abs2Rel}
Let $\FF$ be a totally ordered division ring, let $n<\omega$.
If a finite lattice $L$ embeds into $\mathbf{K}(\FF^n)$, then
it embeds into $\Co(\FF^n,\Omega)$ for all large enough finite
$\Omega\subset\FF^n$.
\end{corollary}

Now let $x$, $a_0$, $a_1$, $b_0$, $b_1$, $c_0$, $c_1$ be variables,
define new terms by
       \begin{align}
       x'&=x\wedge(a_0\vee a_1)\wedge(b_0\vee b_1)\wedge(c_0\vee
c_1),\label{Eq:x'}\\
       a_{i,j,k}&=a_{1-i}\vee\pI{(a_i\vee x')\wedge(b_j\vee
c_k)},\label{Eq:aijk}\\
       b_{i,j,k}&=b_{1-j}\vee\pI{(b_j\vee x')\wedge(a_i\vee
c_k)},\label{Eq:bijk}
       \end{align}
and consider the following lattice-theoretical identity:
       \begin{equation}\label{Eq:StarId}
       x'=\bigvee_{i,\,j,\,k<2}\pI{(x'\wedge a_{i,j,k})\vee(x'\wedge
b_{i,j,k})}.
       \end{equation}

\begin{lemma}\label{L:NewIdDim2}
The lattice $\Co(\FF^2)$ satisfies the identity~\eqref{Eq:StarId}, for any
totally ordered division ring $\FF$.
\end{lemma}

\begin{proof}[Outline of proof]
Let $X$, $A_0$, $A_1$, $B_0$, $B_1$, $C_0$, $C_1$ in $\Co(\FF^2)$, let $X'$,
$A_{i,j,k}$, $B_{i,j,k}$, for $i$, $j$, $k<2$, be formed from these parameters
as in~\eqref{Eq:x'},~\eqref{Eq:aijk}, and~\eqref{Eq:bijk}. Denote by $Y$
the right hand side of~\eqref{Eq:StarId} formed with these parameters. As
it is obvious that $Y$ is contained in $X'$, it suffices to prove that $X'$
is contained in
$Y$. Let $x\in X'$. If $x\in A_i\cup B_i\cup C_i$, for some $i<2$, then
$x\in Y$; thus suppose that $x\notin A_i\cup B_i\cup C_i$, for all $i<2$.
Since $x\in A_0\vee A_1$, there are $a_i\in A_i$, for $i<2$, such that
$x\in[a_0,a_1]$. Similarly, there are $b_i\in B_i$ and $c_i\in C_i$, for $i<2$,
such that $x\in[b_0,b_1]\cap[c_0,c_1]$. Observe that
$x\notin\set{a_i,b_i,c_i}$, for all $i<2$.

Let $\ell$ be the affine line containing $\set{c_0,c_1}$,
and let $i$, $j<2$ such that $a_i$ and $b_j$
are on one side of $\ell$ while
$a_{1-i}$ and $b_{1-j}$ are on the other side.
Take $x$ as origin of the affine plane, and pick any affine line $\ell'$
such that $x\in\ell'$ and either both $a_i$ and $b_j$ are on $\ell'$ (if $x$,
$a_i$, $b_j$ are collinear) or $a_i$ and $b_j$ are on opposite sides of
$\ell'$ (otherwise). Take $(\ell,\ell')$ as a coordinate system in
which $a_i$ and $b_j$ have $\ell'$-coordinates at least~$0$ while $a_{1-i}$ and
$b_{1-j}$ have $\ell'$-coordinates at most~$0$.
Expressing $x$, $c_0$, $c_1$, $a_i$, $b_j$ in this coordinate system yields,
up to possible permutation of $(a_0,a_1)$ and $(b_0,b_1)$, an integer $k<2$ and
elements $\alpha$, $\beta$, $\alpha'$, $\beta'$ of $\FF^+$ and  $\gamma_0$,
$\gamma_1\in\FF^{++}$ such that $\alpha'\leq\beta'$ and
       \begin{align*}
       a_i&=(-\alpha,\alpha'),&b_j=(\beta,\beta'),\\
       c_{1-k}&=(-\gamma_0,0),&c_k=(\gamma_1,0),\\
       x&=(0,0).
       \end{align*}
A careful inspection of every case yields that $[a_i,c_k]\cap[x,b_j]$ is
always nonempty. If $z$ denotes any element of this set, then $x$ belongs
to $[b_{1-j},z]$, thus to $B_{i,j,k}$, thus to~$Y$.
\end{proof}

\begin{lemma}\label{L:SevenPt}
There exists a seven-element subset $\Omega$ of $\QQ^2$ such that
$\Co(\QQ^2,\Omega)$ does not satisfy the identity~\eqref{Eq:StarId}.
\end{lemma}

\begin{proof}
Put $\Omega=\set{\ad_0,\ad_1,\bd_0,\bd_1,\cd_0,\cd_1,\xd}$, where
       \begin{align*}
       \ad_0&=(-2,0),&\ad_1&=(2,0),\\
       \bd_0&=(-1,1),&\bd_1&=(1,-1),\\
       \cd_0&=(1,1),&\cd_1&=(-1,-1),\\
       \xd&=(0,0).
       \end{align*}
Put $x=\set{\xd}$, $a_i=\set{\ad_i}$, $b_i=\set{\bd_i}$, $c_i=\set{\cd_i}$, for
all $i<2$. Then it is straightforward to compute that with those parameters,
the right hand side of~\eqref{Eq:StarId}, calculated in $\Co(\QQ^2,\Omega)$, is
empty, while the right hand side is $x$. Hence $\Co(\QQ^2,\Omega)$ does not
satisfy~\eqref{Eq:StarId}.
\end{proof}

\begin{corollary}\label{C:SevenPt}
Let $\Omega$ be the seven-element set of Lemma~\textup{\ref{L:SevenPt}}. Then
$\Co(\QQ^2,\Omega)$ cannot be embedded into $\Co(\FF^2)$, for any totally
ordered division ring $\FF$.
\end{corollary}

Other phenomena may happen. For example, if $C$ is a square of
$\QQ^2$ (e.g., see
Example~\ref{Ex:PlnonPl}) and $C'=C\cup\set{c}$ where $c$ is the center of $C$,
then $\Co(\QQ^2,C)\cong\two^4$ has a lattice embedding into $\mathbf{K}(\QQ^2)$
(send every $a\in C$ to the segment $[a,c]$), but it has no
zero-preserving such
embedding. On the other hand, $C'$ is a plenary subset of $\QQ^2$ (see
Definition~\ref{D:PlSub}), thus $\Co(\QQ^2,C')$ has a plenary zero-preserving
lattice embedding into $\Co(\QQ^2)$. Observe that $\Co(\QQ^2,C)$ is a
homomorphic
image of $\Co(\QQ^2,C')$.

\section{Open problems}\label{S:Pbs}

In view of Theorem~\ref{T:ArbLatt}, it is natural to ask whether any
finite lattice embeds
into $\mathbf{K}(\QQ^n)$, for some natural number $n$. However, the latter
lattice is known to be \jsd.

\begin{problem}\label{Pb:FinSD+}
Is it the case that every finite \jsd\ lattice can be embedded into
$\mathbf{K}(\QQ^n)$, for some natural number $n$?
\end{problem}

By Theorem~\ref{T:FinLB}, Problem~\ref{Pb:FinSD+} can be answered positively
for finite lower bounded lattices.

Define \emph{semi-algebraic} convex subsets of $\QQ^n$ to be the solution sets
of finite systems of linear inequalities (allowing both $\leq$ and $<$), that
is, the finite intersection of either open or closed affine half-spaces of
$\QQ^n$.

\begin{problem}\label{Pb:FinLatt}
Can every finite lattice be embedded into the lattice $\ol{\mathbf{K}}(\QQ^n)$
of \emph{bounded semi-algebraic} convex subsets of $\QQ^n$, for some natural
number $n$?
\end{problem}

It is well-known that for every Hausdorff locally convex topological vector
space~$V$ over $\RR$, the lattice $\mathbf{CB}(V)$ of all \emph{convex bodies}
of $V$, that is, compact convex subsets of $V$, is \jsd. The proof is
analogous to the one of~\cite[Theorem~15]{BB}.
\begin{problem}\label{Pb:SD+}
Is it the case that every \jsd\ lattice can be embedded into $\mathbf{CB}(V)$,
for some Hausdorff locally convex topological vector space $V$ over $\RR$?
\end{problem}

\begin{problem}\label{Pb:ClosedConv}
Is it the case that every lattice can be embedded into the lattice of all
\emph{bounded closed convex} subsets of some real Banach space?
\end{problem}

Our next problem asks about dependence from the division ring $\FF$. It follows
from Theorem~\ref{T:ArbLatt} that for a totally ordered division ring
$\FF$, the
universal theory, in the language $(\vee,\wedge)$, of $\Co(\FF^{(I)})$, for
infinite $I$, is the universal theory of all lattices. This leaves open the
problem in finite dimension.

\begin{problem}\label{Pb:DepFF}
For a natural number $n$ and a totally ordered division ring $\FF$, do the
lattices $\Co(\QQ^n)$ and $\Co(\FF^n)$ have the same universal theory?
\end{problem}

It follows from~\cite{SeWe3}
that the answer to
Problem~\ref{Pb:DepFF} is positive for $n=1$. Also observe that
$\Co(\QQ^2)$ and $\Co(\RR^2)$ do not have the same first-order theory, see
B. Gr\"unbaum~\cite[Example~5.5.3]{Grun}.

\section*{Acknowledgment}
The results of the paper were obtained while both authors were guests
of the Department of Algebra at Charles University in Prague, in June
2002. The excellent conditions provided there cannot be forgotten.
Special thanks are due to Ji\v{r}\'\i\ T\r{u}ma and V\'aclav Slav\'\i k.

\end{document}